\documentclass[a4paper,abstract=true]{scrartcl}
\usepackage[utf8]{inputenc}

\usepackage[english]{babel} 
\usepackage{amsmath, amsthm, amssymb} 
\usepackage{csquotes, geometry, enumitem} 
\usepackage{xspace, xparse, graphicx} 
\usepackage{tabu} 
\usepackage{tikz-cd}
\usepackage{stmaryrd}
\usepackage{url}
\usepackage{appendix}

\usepackage[hidelinks]{hyperref}

\theoremstyle{definition}
\newtheorem{defn}{Definition}[section]
\newtheorem{rmk}[defn]{Remark}
\newtheorem{ex}[defn]{Example} 

\theoremstyle{plain}
\newtheorem{thm}[defn]{Theorem}
\newtheorem*{thm*}{Theorem}
\newtheorem{lem}[defn]{Lemma}
\newtheorem{prop}[defn]{Proposition}
\newtheorem{cor}[defn]{Corollary}

\newlist{HS}{enumerate}{1} 
\setlist[HS]{label=HS(\arabic*)}

\newlist{SB}{enumerate}{1}
\setlist[SB]{label=SB(\arabic*)}

\NewDocumentCommand\party{mmO{}}{\frac{\partial^{#3} #1}{\partial^{#3} #2}}
\NewDocumentCommand\dx{}{\operatorname{d\!}}

\NewDocumentCommand\Id{}{\operatorname{Id}}
 
\NewDocumentCommand\del{}{\partial} 
\NewDocumentCommand\delb{}{\bar{\partial}}

\NewDocumentCommand\Znum{}{\mathbb{Z}}
\NewDocumentCommand\Rnum{}{\mathbb{R}}
\NewDocumentCommand\Cnum{}{\mathbb{C}}

\NewDocumentCommand\Symp{}{\operatorname{Symp}}
\NewDocumentCommand\Diff{}{\operatorname{Diff}}

\newcommand{\cM}{\mathcal{M}}
\newcommand{\cL}{\mathcal{L}}
\newcommand{\cP}{\mathcal{P}}

\newcommand{\ZZ}{\ensuremath{\mathbb Z}}
\newcommand{\CC}{\ensuremath{\mathbb C}}
\newcommand{\RR}{\ensuremath{\mathbb R}}

\newcommand{\QQQ}{\mathbf{Q}}
\newcommand{\Qw}{\mathcal{Q}_{[\omega]}}
\newcommand{\HtR}{(H^{2,0}\oplus H^{0,2})_{\Rnum}(M)}

\title{Moduli spaces of spacefilling branes\\in symplectic 4-manifolds}
\author{Charlotte Kirchhoff-Lukat\thanks{KU Leuven, Department of Mathematics, Celestijnenlaan 200B box 2400, BE-3001 Leuven, Belgium.
Email: 
\texttt{charlotte.kirchhofflukat@kuleuven.be}}, Marco Zambon\thanks{KU Leuven, Department of Mathematics, Celestijnenlaan 200B box 2400, BE-3001 Leuven, Belgium.
Email: 
\texttt{marco.zambon@kuleuven.be}}}
\date{}

\begin{document}

\maketitle

\vspace{-1cm}
\begin{abstract} 
On a symplectic manifold $(M,\omega)$, a spacefilling brane structure is a closed 2-form $F$ which determines a complex structure, with respect to  which  $F+i\omega$ is holomorphic symplectic. For 
holomorphic symplectic compact Kähler 4-manifolds,  
we show that the moduli space of spacefilling branes is smooth, and determine its dimension. The proof relies on the local Torelli theorem for K3 surfaces and tori. 
\end{abstract}

\smallskip
\begin{center}
    \emph{\large{Dedicated to the memory of our student Gilles Castel}}
\end{center}
\bigskip
%\bigskip

\setcounter{tocdepth}{2}%doesn't disply subsubsections
\tableofcontents

\section{Introduction} 
This paper is concerned with a specific geometric structure defined on a symplectic manifold $(M,\omega)$, namely \emph{spacefilling brane} structures. They are provided by closed $2$-forms $F$  such that the endomorphism  $$I:=\omega^{-1} \circ F\colon TM\to TM$$ squares to $-Id$. In this case, $I$ is an integrable complex structure, and
$F+i\omega$ is a \emph{holomorphic symplectic} form (in particular, $\omega$ is of type\footnote{This is in  contrast to Kähler forms, which are real of type $(1,1)$.}  $(2,0)+(0,2)$ w.r.t. $I$, and $\dim_{\RR}(M)$ is a multiple of $4$). 
Conversely, on the symplectic manifold $(M,\omega)$,  a complex structure $I$ can be used to construct a holomorphic symplectic form -- and thus a spacefilling brane -- if and only if $F(\cdot,\cdot):= \omega(I\cdot,\cdot)$ is a 2-form, i.e. \emph{antisymmetric}\footnote{This is in contrast with K\"ahler geometry, where $ \omega(I\cdot,\cdot)$ is symmetric and positive definite.}. The space of spacefilling branes for a fixed symplectic form $\omega$ can thus be viewed as the space of precisely those holomorphic symplectic forms with fixed imaginary part $\omega$.

The notion of \emph{brane}  comes from string theory, where they are, very loosely speaking, submanifolds of spacetime that strings can end on, equipped with their own structure. The A-model of string theory is defined in terms of the symplectic structure of the spacetime manifold, while the B-model is defined in terms of its complex structure. In both cases, branes are natural submanifolds for these structures: Lagrangian -- and, as first observed by \cite{KO} -- certain higher-dimensional coisotropic submanifolds, for the A-model; complex submanifolds (equipped with suitable local systems) for the B-model. 
Generalized complex geometry offers a unified description for both of these settings, as well as a general definition of what a brane should be  \cite{GualtieriThesis, KaLi05}. Branes, including the spacefilling ones,  on general generalized complex manifolds have garnered substantial research interest in recent years, see for example \cite{Gu10, BG21}.

\bigskip

We will restrict ourselves to the symplectic setting, and no knowledge of generalized complex geometry is assumed. Higher-dimensional coisotropic branes in symplectic manifolds are much less studied than the lowest-dimensional branes, Lagrangian submanifolds. We are interested in the deformation theory of coisotropic branes: one reason is that it provides new examples of branes, another is that it is  closely related to the deformation theory of coisotropic submanifolds, which is rich and highly non-linear \cite{OP}.

In this paper, we focus on spacefilling branes on certain  4-dimensional symplectic manifolds. We go beyond their deformation theory and  are able to give a \emph{global} description of the moduli space. In the companion paper \cite{PrequantisedSpacefillingBranes}
we will consider prequantised branes, study their existence, and compare their deformations to those appearing in this paper.

\bigskip
\noindent{\textbf{Main results}}

In this note we  chiefly  focus on compact  $4$-manifolds   
admitting 
a complex structure which has trivial canonical bundle and which is 
K\"ahler.
 There are two such classes of complex manifolds: Complex tori $\Cnum^2/\Lambda$ (where $\Lambda\cong \Znum^4$), and K3 surfaces, the simply connected compact complex Kähler surfaces with trivial canonical bundle. For each of these two classes, the underlying smooth manifolds are all diffeomorphic, being the 4-torus $T^4$ and the K3 manifold  respectively, but the complex structures are in general non-equivalent.  

Let $M$  be either $T^4$ or the K3 manifold, and fix a symplectic form $\omega$ admitting a spacefilling brane.
 Our main result is a statement about the moduli space $\cM_{\omega}$, the quotient of the space of spacefilling branes  on $(M,\omega)$ by the group $\Symp_*(M,\omega)$ of symplectomorphisms inducing the identity in cohomology.
To state the result, consider the map
\begin{align}
\label{eq:globtorintro}
\Phi \colon \cM_{\omega}&\to \mathcal{Q}_{[\omega]} \\
 F' \text{ mod } \Symp_*(M,\omega) &\mapsto [F'], \nonumber
\end{align}
 where  $\mathcal{Q}_{[\omega]}$ is 
 the codimension two submanifold of  $H^2(M,\RR)$ cut out by the conditions $[F']\wedge [\omega]= 0$ and $[F'] \wedge [F']=[\omega] \wedge [\omega]$.
 
 We paraphrase Theorem \ref{prop:CMHR} and Corollary \ref{lem:compact}.

\begin{thm*} 
{Let $M$ be the K3 manifold or the 4-torus, endowed with any symplectic form $\omega$ admitting a spacefilling
brane structure.}

The moduli space of spacefilling brane structures $\cM_{\omega}$  is a non-Hausdorff smooth manifold, and non-compact. 
{If $M$ is the K3 manifold, the moduli space $\cM_{\omega}$  is of dimension $20$; if $M$ is  the  4-torus, the moduli space is of dimension $4$.}
%For $M$ a K3 surface (respectively the 4-torus), endowed with any symplectic form $\omega$, the moduli space $\cM_{\omega}$  is of dimension $20$ (respectively $4$). 

The map $\Phi$ of \eqref{eq:globtorintro} is smooth, and is a local diffeomorphism.

\end{thm*}

In particular, given a spacefilling brane $F$, there is a neighborhood of $[F]$ in $\mathcal{Q}_{[\omega]}$ such that 
any point of the neighborhood corresponds to a class of spacefilling branes $F'$, which are not related to $F$ by any element of $\Symp_*(M,\omega)$.

\bigskip
The above result relies in an essential way on a version of the Local Torelli Theorem, {as we now outline.} 

First  we encode spacefilling branes as {those} complex structures $I$ such that $\omega(I\cdot,\cdot)$ is skew-symmetric, {via the assignment $F \mapsto I:=\omega^{-1}\circ F$, see Proposition \ref{prop:charbrane} ii)}.

On a compact  $4$-manifold   admitting 
%a brane structure $F$ for which the corresponding complex structure is K\"ahler. 
a complex structure which has trivial canonical bundle, there is a well-known way to encode any complex structure  by a line of closed complex 2-forms, {which we recall in Proposition \ref{prop:corr1}.
Complex structures satisying the above skew-symmetry property correspond  to complex 2-forms with imaginary part $\omega$.}

The Local Torelli Theorem as stated in \cite{HuyTalk}, which holds under the
K\"ahler condition and we  recall in Theorem \ref{thm:loctor},  is a statement about the induced map from the moduli space of complex structures to the projectivization of $H^2(M,\CC)$. A suitable restriction of this map yields the map $\Phi$    in \eqref{eq:globtorintro}.

 \bigskip
\noindent{\textbf{Notation:}} 
Given a 2-form $\sigma$, we denote by the same symbol the induced bundle map $TM\to T^*M, \;v\mapsto \iota_v\sigma.$

\bigskip
\noindent{\textbf{Acknowledgements:}} 
We thank Marco Gualtieri  for crucial input about the local Torelli theorem he gave us at the beginning of this project, and Luka Zwaan for pointing us to the existence of smooth universal deformations for the proof of Lemma \ref{lem:3equiv}. We also thank Ignacio Barros for interesting discussions, {Jim Stasheff for feedback, and the referee for valuable suggestions that improved the manuscript.}

C.K-L. acknowledges support by the FWO Junior Postdoctoral Fellowship 1263121N and Marie Skłodowska-Curie Fellowship 887857 under the European Union's Horizon 2020 research framework. 
M.Z. acknowledges partial support by
the FWO and FNRS under EOS projects G0H4518N and G0I2222N, by FWO project G0B3523N, and by Methusalem grant METH/21/03 - long term structural funding of the Flemish Government (Belgium).

\section{{Spacefilling} branes in symplectic manifolds} 

In this section we review spacefilling branes in symplectic manifolds, {establishing} several {equivalent} characterizations. The {general notion} of brane arises naturally in the context of {both string theory and} generalized complex geometry \cite[\S 3]{KaLi05}.
%, and we  recall it later in Definition \ref{defn:branegen}. 
Specialized to the symplectic context, it reads as follows:

\begin{defn} \label{def:cobrane} (\cite{GualtieriThesis}, Example 7.8)
A  \emph{coisotropic A-brane} in the symplectic manifold $(M,\omega)$ is a pair $(Y,F)$ of a submanifold {$i_Y\colon Y \hookrightarrow  M$} and a closed 2-form $F\in \Omega^2(Y,\RR)$ s.t. 
\begin{enumerate}[label=(\roman*)]
\item $Y$ is a coisotropic submanifold, 
\item $\ker(F)=\ker(i_Y^*\omega)=: T^{\omega}_Y$, 
\item the endomorphism $I:= \tilde{\omega}^{-1} \circ \tilde{F}$  squares to $-Id$, where  we denote by $\tilde{F},\tilde{\omega}$ the descended forms to $TY/T_Y^{\omega}$.
\end{enumerate}
\end{defn} 

\begin{rmk}
It follows that $I$ is a transverse complex structure, i.e. on any submanifold 
transverse to $T^{\omega}_Y$ it is integrable. The fact that    $I$ is a transverse complex structure is shown in \cite[\S 4]{KO}, and we will provide a direct computational argument below.

Note that in the case where the quotient of $Y$ by its characteristic foliation (the foliation to which $T_Y^{\omega}$ is tangent) is a smooth manifold, {this leaf space} is holomorphic symplectic: 
By construction, $\tilde{F} + i\tilde{\omega}$ is a closed holomorphic 2-form with respect to the complex structure $I=\tilde{\omega}^{-1}\circ \tilde{F}$, and the integrability of the latter is equivalent to the closedness of $\tilde{F},\tilde{\omega}$.
(Here, by a slight abuse of notation, we denote by $\tilde{F},\tilde{\omega}$ the induced 2-forms on the leaf space.)
Recall that a closed holomorphic 2-form is \emph{holomorphic symplectic} when its imaginary part is non-degenerate.

\end{rmk}

\begin{rmk}
    The systematic study of coisotropic A-branes in symplectic manifolds is originally due to Kapustin-Orlov in $\cite{KO}$, who establish that there are higher-dimensional branes supported on some coisotropic submanifolds.  {The description as in Definition \ref{def:cobrane} can be found in \cite[§4]{KO}.}
\end{rmk}

In this article, we focus on the case where a brane structure is supported on an entire symplectic manifold, i.e. the manifold is equal to the leaf space. The above definition then specializes to the following: 

\begin{defn} 
A \emph{spacefilling brane structure} on a symplectic manifold $(M,\omega)$ is a closed 2-form $F\in\Omega^2_{\operatorname{cl}}(M,\RR)$ such that $\omega^{-1} \circ F$ squares to $-Id$.
\end{defn} 

Here, by a slight abuse of notation, we denote by $\omega$ also the corresponding bundle map $TM\to T^*M$, and similarly for $F$.

In this case, $\omega^{-1} \circ F$ defines a complex structure on $M$, with respect to which $F+i\omega$ is  holomorphic symplectic. We will show this in Proposition \ref{prop:charbrane}; to this aim, we state and prove the following known result  {(see \cite[\S 4]{KO} for one direction of this equivalence)}: 

\begin{lem}\label{prop:inte}
    Let $(M,\omega)$ be a symplectic manifold and $I$ an almost complex structure on $M$ s.t. $F:=\omega \circ I$ is a 2-form, i.e. antisymmetric. 

    Then {the following holds:} $I$ is a complex structure, i.e. integrable, if and only if $\dx F=0$. 
\end{lem}
%Note that $M$ in this Proposition should be viewed as the leaf space of a coisotropic A-brane.
\begin{proof}

Denote by $N_I$ the Nijenhuis tensor of the almost complex structure $I$. By the Newlander-Nirenberg theorem, $I$ is integrable iff $N_I=0$. 

We begin by computing: 
For $X,Y\in\Gamma(TM)$: 
\begin{align*}
    N_I(X,Y) 
    &= [X,Y] + I \left([IX,Y] + [X,IY]\right) - [IX,IY] \\ 
    &= [X,Y] - I\left(\cL_Y (IX)+\cL_{IY}X\right) + \cL_{IY}(IX) \\ 
    &=[X,Y] - I\left( (\cL_YI)(X) + I \cL_Y X + \cL_{IY} X\right) + (\cL_{IY}I)(X) + I \cL_{IY}X \\ 
    &= [X,Y] - I(\cL_YI)(X) + [Y,X]+(\cL_{IY}I)(X) \\ 
    &= \left(\cL_{IY}I-I\cL_Y I\right)(X)
\end{align*}
Thus  we have the following characterization\footnote{A similar characterization, in which $\cL$ is replaced by the Levi Civita connection of a Riemannian metric preserved by $I$, can be found e.g. in  \cite[Lemma 5.4]{Moroianu_2007}.} of the vanishing of $N_I$:

\begin{equation}\label{eq:NIzero}
    N_I =0\ \Leftrightarrow\ \cL_{IY} I - I\cL_Y I=0\quad \forall Y \in \Gamma(TM).
    \end{equation} 

Now we compute for all vector fields $Y\in \Gamma(TM)$:
\begin{align*}
    \cL_Y I &= \cL_Y(\omega^{-1}\circ F) \\ 
    &= (\cL_Y \omega^{-1})\circ F + \omega^{-1} \circ \cL_Y F \\ 
    &= -\omega^{-1}(\cL_Y\omega) \circ (\omega^{-1}\circ F) + \omega^{-1} \circ  \cL_Y F \\ 
    &= \omega^{-1} \left(-(\dx i_Y\omega)\circ I + \dx i_Y F + i_Y \dx F\right). 
\end{align*}
In the second equality we applied the product rule for the Lie derivative of tensors. 

Thus the right hand side of eq. \eqref{eq:NIzero} can be written as: 
\begin{align*}
    \cL_{IY}I - I\cL_Y I =& \cL_{IY} I + (\cL_Y I)I \\ 
    =&\omega^{-1}\left(-(\dx i_{IY}\omega) \circ I + \dx i_{IY}F + i_{IY}\dx F \right. \\ 
    &\quad\quad\left.- \dx i_Y \omega \circ I^2 + (\dx i_Y F)\circ I + (i_Y \dx F)\circ I \right) \\ 
    %&= \omega^{-1} \left(-(\dx i_Y F)\circ I - \dx i_{IY}\omega + i_{IY}\dx F \right. \\ 
    %&\left.+\dx i_Y \omega + (\dx i_Y F)\circ I + (i_Y \dx F)\circ I \right) \\ 
    =& \omega^{-1} \left(i_{IY}\dx F + (i_Y \dx F) \circ I\right),
\end{align*}
using $\omega\circ I =F$ in the last equality.
Summarizing, we have shown that 
\begin{equation*} 
    N_I =0\ \Leftrightarrow\ i_{IY}\dx F + (i_Y \dx F) \circ I=0\quad \forall Y \in \Gamma(TM).
    \end{equation*} 

Consequently, if $\dx F=0$, it follows that $N_I=0$. 

Conversely, if $N_I=0$, 
we find that 
\[ 
(\dx F)(IX, Y,Z)+(\dx F)(X, IY, Z) = 0\quad \forall X,Y,Z \in \Gamma(TM). 
\]
In particular, if we take $X=\party{}{z^i}, Y = \party{}{z^j}$ for local holomorphic coordinates, and $i\neq j$, we find 
\[
    (\dx F)\left(\party{}{z^i},\party{}{z^j}, Z\right)=0\quad \forall Z\in\Gamma(TM). 
\] 
Similarly, taking $X=\party{}{\bar{z}^i}, Y = \party{}{\bar{z}^j}$: 
\[
    (\dx F)\left(\party{}{\bar{z}^i},\party{}{\bar{z}^j}, Z\right)=0\quad \forall Z\in \Gamma(TM).
\] 
This implies that $\dx F=0$.\end{proof}

\begin{prop}\label{prop:charbrane}

Let $(M,\omega)$ be a symplectic manifold. The following structures are in bijection with each other:
\begin{itemize}
    \item[i)] A spacefilling brane structure $F$,
    \item[ii)] A complex structure $I$ such that $\omega \circ I$ is skew-symmetric
    \item[iii)] A complex structure $I$ and a 2-form $F$ such that $F+i\omega$ is holomorphic symplectic w.r.t. $I$.
\end{itemize}

\end{prop}

\begin{proof}
The bijection between $i)$ and $ii)$ is given by $$F\mapsto I:=\omega^{-1} \circ F \quad \text{ and }\quad I\mapsto F:=\omega \circ I.$$ Given $F$, we know that $I$ is almost complex, and  the integrability of $I$ is shown in Lemma \ref{prop:inte}.
Given $I$, we know that $F$ is a 2-form, and its closedness follows from Lemma \ref{prop:inte}.

It is easy to see that $iii)$ implies the other two items: $F+i\omega$ being of type $(2,0)$ implies that $F=\omega \circ I$, and the closeness of $F+i\omega$ implies the one of $F$. Given $i)$ and $ii)$, the equation $I=\omega^{-1} \circ F$ implies that $F+i\omega$ is of type $(2,0)$, and the closeness of $F$ and $\omega$ implies that $F+i\omega$ is closed, in particular $\bar{\partial}$-closed, thus holomorphic.   
\end{proof}

\begin{rmk}\label{rem:holDarboux}
The complex structure $I$ in \emph{ii)} induces the same orientation on $M$ as the symplectic form $\omega$. This can be seen writing the holomorphic symplectic form $F+i\omega$ in holomorphic Darboux coordinates (see e.g.   \cite{MoserHolSympl}); for $\dim(M)=4$, we display explicit formulae in Remark \ref{ex:torus}.
\end{rmk}
 
 \section{On complex 4-manifolds}
\label{sec:complex4man}

In this section we review some facts about complex structures on compact $4$-manifolds. 

A feature of dimension $4$ is that
the wedge product on the space of 2-forms induces a pointwise inner product. On 4-manifolds which admit nowhere vanishing holomorphic 2-forms, this allows for a convenient characterisation of complex structures in terms of complex 2-forms
(see \S \ref{subsec:4man}). Furthermore, for Kähler 4-manifolds, 
the moduli space of complex structures can -- at least locally -- be described in terms of de Rham cohomology (see \S \ref{subsec:localtorelli}).

\subsection{Complex 4-Manifolds with a holomorphic symplectic form}\label{subsec:4man}

We start by recalling a well-known fact, which will be needed to state Proposition \ref{prop:corr1}.

\begin{prop}
\label{prop:HS}
(See e.g. Proposition 1.2 in \cite{HuyLecturesK3}.) Let $M$ be an oriented 4-manifold.
Given a complex 2-form $\Omega\in \Omega^2(M,\Cnum)$, the conditions that $\Omega$ be holomorphic symplectic with respect to \emph{some} complex structure on $M$ defining the given orientation are: 
\begin{HS}
\item $\Omega\wedge \bar{\Omega}> 0$ everywhere on $M$ (non-degeneracy)
\item $\Omega\wedge\Omega= 0$ (holomorphicity; a real 4-manifold cannot admit holomorphic 4-forms), 
\item $\dx \Omega =0$ (closedness). 
\end{HS}
\end{prop}

Given a form $\Omega$ that satisfies these properties, the complex structure with respect to which it is holomorphic is uniquely determined. Indeed, it is determined by the kernel of $\Omega$, which is integrable because we assume $\Omega$ to be closed: 
\begin{equation}\label{eq:T10}
    T^{0,1}M:= \{X\in T_{\Cnum}M \vert i_X \Omega=0\}
\end{equation} 
 Fixing the antiholomorphic tangent bundle is equivalent to fixing a complex structure.
Note that the orientation determined by the volume form $\Omega \wedge \bar{\Omega}$ agrees with the one induced by the corresponding complex structure.

\begin{rmk}\label{rmk:ortho}
On an orientable 4-manifold $M$, a choice of volume form yields an identification 
$\wedge^4T^*_pM \cong \Rnum$ at every point $p$. The wedge product on $\wedge^2 T^*M$ thus defines an inner product  on every $\wedge^2 T_p^*M$. Here by \emph{inner product} we mean a  non-degenerate symmetric bilinear form,  {possibly of indefinite signature.}  

For later use in \S \ref{sec:defdiff}, 
we introduce the following notation:
Fixing a complex structure on $M$, denote by $\Omega^{1,1}_{\RR}(M)$ or $\Omega^{1,1}_{\RR}$ the real elements of $\Omega^{1,1}(M,\Cnum)$, i.e. those that remain invariant under complex conjugation. Similarly let $(\Omega^{2,0}\oplus \Omega^{0,2})_{\RR}(M)$ or simply $(\Omega^{2,0}\oplus \Omega^{0,2})_{\RR}$ denote the real elements of $\Omega^{2,0}(M,\Cnum)\oplus \Omega^{0,2}(M,\Cnum)$. Call $(T^*)^{(2,0)+(0,2)}_{\RR}$ the real vector subbundle of $\wedge^2 T^*M$ whose sections are $(\Omega^{2,0}\oplus \Omega^{0,2})_{\RR}$, and define similarly $(T ^*)^{(1,1)}_{\RR}$.

 The direct sum
$$\wedge^2 T^*M= (T ^*)^{(1,1)}_{\RR} \oplus (T^*)^{(2,0)+(0,2)}_{\RR}$$ is an orthogonal decomposition   with respect to the pointwise inner product induced by $\wedge$, since there are no $(3,1)$-forms or $(1,3)$-forms on a 4-manifold.
\end{rmk}

\begin{rmk}\label{rmk:K}  
As a preparation for the proof of Proposition \ref{prop:corr1}, we 
recall that on an 4-manifold $M$, any complex structure $I$  defines a holomorphic line bundle $K:=\wedge^2T^*_{1,0}M\subset \wedge^2 T^*_{\Cnum}M$ of 2-forms which are of type $(2,0)$ with respect to $I$, the \emph{canonical line bundle}. The canonical line bundle $K$ is trivial as a holomorphic line bundle if{f} it admits a nowhere-vanishing holomorphic section.  
Such a nowhere-vanishing holomorphic section is necessarily  a holomorphic symplectic form, since in complex  local coordinates it is of the form $f dz_1\wedge dz_2$ with $f$ a nowhere-vanishing holomorphic function. 
Conversely, the existence of a holomorphic symplectic form implies that the canonical line bundle of the underlying complex structure is trivial.

Note: If $M$ is compact, this holomorphic section is unique up to multiplication by a constant $c\in \Cnum\setminus \{0\}$,  
since the only holomorphic functions on a compact complex manifold are constant.
\end{rmk}

 From now on, assume that all complex structures induce the chosen global orientation. 

We recall a well-known correspondence, which we will use in \S \ref{subsec:localtorelli}.   {In order to be self-contained, we provide a proof}. As mentioned in \cite[Remark 0.1]{HuyTalk}, for K3 surfaces this description appears to be due to Andreotti.

 \begin{prop} \label{prop:corr1}
 On an oriented compact 4-manifold $M$ which admits a complex structure and a holomorphic symplectic form, there is a one-to-one correspondence: 
\[
\left\{
\begin{matrix}
\text{Complex structures}\\
\text{defining the given orientation}
\end{matrix}
\right\}
\stackrel{1:1}{\longleftrightarrow} 
\left\{
\begin{matrix}
\text{Complex lines of 2-forms }\\ \Cnum\Omega \subset \Omega^2(M,\Cnum) \text{ satisfying HS(1)-HS(3)}
\end{matrix}
\right\} 
\]
 \end{prop} 

\begin{proof}
The map \enquote{$\leftarrow$}  {associating a complex structure to each complex 2-form satisfying HS(1)-HS(3)} is described in the text just after Proposition  \ref{prop:HS}. 

The map \enquote{$\rightarrow$} assigns to a complex structure the line of 2-forms which are  holomorphic symplectic with respect to it. 
To explain this, we use the fact -- proven below -- that every complex structure on the smooth compact manifold $M$ has trivial canonical bundle.  By Remark \ref{rmk:K}, the canonical bundle thus has a nowhere-vanishing holomorphic section $\Omega$, which is necessarily a holomorphic symplectic form, and all holomorphic sections are multiples of $\Omega$ by a constant $c\in \CC\setminus{\{0\}}$. If the complex structure defines the given orientation on $M$, then  {so does} $\Omega$. 

{The Kodaira classification tells us that the canonical bundle of \emph{every} complex structure on the smooth compact manifold $M$ is indeed trivial:} 
%\emph{every} complex structure on the smooth compact manifold $M$ has trivial canonical bundle. We can argue this using the Kodaira classification: 
The only compact complex surfaces with trivial canonical bundle are K3 surfaces, complex tori $\Cnum^2/\Lambda$ (both always Kähler) and primary Kodaira surfaces (never Kähler). These are fully characterised by their Kodaira dimension and Hodge diamonds, which in this case are invariants of their underlying smooth manifolds only.  
\end{proof}

\begin{rmk} 
From the Kodaira classification of complex surfaces we see that for a compact complex 4-manifold the condition that the canonical bundle be trivial, and thus a holomorphic symplectic form exist, is very restrictive and only a few classes of examples exist. 

In the following subsection we will examine the Kähler case, where the two classes of examples are complex tori and K3 surfaces; see Remark \ref{rem:holsymK4}. 
The non-Kähler case is less well studied. An examination of compact holomorphic symplectic manifolds that are not Kähler, \emph{primary Kodaira surfaces}, can be found in \cite{BKKY22}. 

% For a compact complex 4-manifold, the condition that the canonical bundle be trivial, and thus a holomorphic symplectic form exist, is very restrictive and only a few classes of examples are known. 

% In the following subsection we will examine the Kähler case, where only 2 classes of examples exist, see Remark \ref{rem:holsymK4}.
%  The non-Kähler case is less well studied. An examination of one class of compact holomorphic symplectic manifolds that are not Kähler, \emph{primary Kodaira surfaces}, can be found in \cite{BKKY22}
\end{rmk}

\subsection{Kähler 4-manifolds and the local Torelli theorem}
\label{subsec:localtorelli}

The main result of this subsection is Theorem \ref{thm:loctor}.
{Now} we focus on compact Kähler 4-manifolds $(M,I)$;
 by this we mean that there exists a real symplectic form $\sigma$ on $M$ which is compatible with the complex structure $I$, in the sense that $\sigma(\cdot,I\cdot)$ is a Riemannian metric. In that case, Hodge theory tells us that the wedge product $\wedge$ defines an inner product -- meaning a non-degenerate symmetric bilinear form, in this case of indefinite signature -- not only pointwise on 2-forms 
 but also on second cohomology $H^2(M,\RR)$. 
 Additionally, the decomposition of its complexification
\begin{equation}\label{eq:Hodge}
H^2(M,\Cnum)= H^{2,0}(M) \oplus H^{1,1}(M) \oplus H^{0,2}(M)
\end{equation}
into Dolbeault cohomology groups is orthogonal with respect to its Hermitian extension.

The inner product space $(H^2(M,\RR),\wedge)$ together with  the above orthogonal decomposition of its complexification form an instance of a \emph{polarised Hodge structure}.
Here and in the following we fix the identification $H^4(M,\RR)\cong \RR$ for which the positive generator of $H^4(M,\Znum)$ gets mapped to $1$ (in other words, the identification given by taking the integral over $M$).

\begin{rmk}\label{rem:holsymK4}
In contrast to the pointwise setting, the signature of this inner product now obviously depends on the global topology of the manifold. There are two classes of compact Kähler 4-manifolds which admit holomorphic symplectic forms: 

 \begin{itemize}
     \item Complex tori $T\cong \Cnum^2/\Lambda$ for a lattice $\Lambda \cong \Znum^4$; in this case, $\dim(H^2(M,\RR))=6$,  
     \item K3 surfaces, which are all  diffeomorphic to the 4-dimensional simply-connected compact solution to the equation $z_0^4+z_1^4+z_2^4+z_3^4=0$ in $\Cnum P^3$ (the \emph{ K3 manifold}),   but can be equipped with many different non-isomorphic complex structures. For  K3 surfaces one has $\dim(H^2(M,\RR))=22$.
 \end{itemize} 
For a complex torus, the inner product on $H^2(M,\RR)$ is of split signature $ (3,3)$, whereas for a K3 surface it has signature $(3,19)$, meaning that it is positive definite on a $3$-dimensional subspace and negative definite on a $19$-dimensional one. 

It turns out that on 4-dimensional tori and on the K3 manifold, \emph{all} complex structures are Kähler: In the case of a complex torus $M\cong \Cnum^2/\Lambda$, the Kähler structure is simply inherited from $\Cnum^2$. In the case of a K3 surface, this is a famous theorem by Todorov and Siu \cite{Tod80,Siu83}. 
\end{rmk} 

Consider first an  oriented compact 4-manifold $M$ which admits a complex structure and a holomorphic symplectic form (equivalently, by Remark \ref{rmk:K}, so that the canonical bundle $K$ is trivial).
{We denote by} $\QQQ$  the subset of $H^2(M,\Cnum)$, cut out by the conditions 
\begin{equation}\label{eq:Q}
    [\Omega'] \wedge [\Omega']=0,\quad [\Omega']\wedge [\bar{\Omega}']> 0. 
\end{equation}
Notice that $\QQQ$ is a complex submanifold\footnote{
This can be seen directly as follows: near any  point of $\QQQ$, we have that $\QQQ$ is the zero level set of the map $[\Omega]\mapsto [\Omega]\wedge [\Omega]\in H^4(M,\Cnum)\cong 
\CC$. Zero is a regular value of this map since $\frac{d}{dt}|_{0}([\Omega]+t[\bar{\Omega}])\wedge ([\Omega]+t[\bar{\Omega}])=2 [\Omega]\wedge [\bar{\Omega}]$ is non-zero.}
of $H^2(M,\Cnum)$. 
Its projectivisation $\mathcal{Q}\subset \mathbb{P}H^2(M,\Cnum)$ is called the \emph{Period domain}, and is a complex submanifold of 
$\mathbb{P}H^2(M,\Cnum$),  because  
the inner product induced on $H^2(M,\Cnum)$ by the wedge product is non-degenerate.

Recall that  Proposition \ref{prop:corr1} provides a correspondence between complex structures $I'$ defining the given orientation and complex lines  $\CC \Omega'$ in $\Omega^2(M,\Cnum)$ satisfying certain conditions.
This correspondence induces a well-defined map 
\begin{equation} \label{eq:globgtor}
\cP\colon
\left\{
\begin{matrix}
\text{Complex structures on $M$}\\
\text{defining the given orientation}
\end{matrix}
\right\}
/ \Diff_*(M)\to \mathcal{Q}:=\mathbb{P}\QQQ, 
\end{equation}
called the \emph{period map},
which to the class of a complex structure $I'$ associates the line $\Cnum[\Omega']$, where $\Omega'$ is a holomorphic symplectic form w.r.t. $I'$. 
Here
\begin{equation}\label{eq:diff*}
    \Diff_*(M):=\left\{\phi \in \Diff(M)\middle\vert \phi^*: H^\bullet (M,\RR)\rightarrow H^\bullet(M,\RR) = \Id_{H^\bullet(M,\RR)}\right\}
\end{equation}
denotes the diffomorphisms inducing the identity in cohomology.

The local Torelli theorem states that the map \eqref{eq:globgtor} is  a local diffeomorphism near the class of $I$, whenever $(M,I)$ is K\"ahler. It gives us a correspondence between small deformations of complex structures and their lines of holomorphic symplectic 2-forms, up to an equivalence relation:

\begin{thm} 
\label{thm:loctor}
\textbf{\emph{(Local Torelli)}}
Let $(M,I)$ be a compact Kähler 4-manifold that admits a holomorphic symplectic structure $\Omega$.

The period map $\cP$ in \eqref{eq:globgtor} locally restricts to diffeomorphism between
\begin{itemize}
    \item small deformations of the complex structure $I$, up to diffeomorphisms in $\Diff_*(M)$ 
    \item small deformations of the complex line $\Cnum[\Omega]$ in $\mathcal{Q}  $.
\end{itemize}
To the class of a complex structure $I'$, the bijection associates the line $\Cnum[\Omega']$, where $\Omega'$ is a holomorphic symplectic form w.r.t. $I'$.
To a complex line $\Cnum[\Omega']$, this bijection associates the class of the complex structure with antiholomorphic tangent bundle determined as in eq. \eqref 
{eq:T10}. 
\end{thm}

For a manifold $M$ as in  Theorem \ref{thm:loctor}, the period map $\cP$ is a \emph{local diffeomorphism}, i.e. it restricts to a diffeomorphism in a neighborhood of any point of the domain. Indeed, every complex structure on $M$ inducing the given orientation is automatically K\"ahler (see Remark \ref{rem:holsymK4}) and admits a holomorphic symplectic form (Proposition \ref{prop:corr1}).
 
\begin{rmk}
    The \emph{local Torelli problem} -- the question of when a period map from equivalence classes of complex structures to cohomology is a local isomorphism -- is a classical problem that has been studied for different classes of manifolds. 
    The provenance of the result for K3 surfaces specifically 
    {(\cite[page 4]{HuyTalk}\cite[Proposition 2.8 in Ch. 6]{HuyLecturesK3})}
    appears to be somewhat unclear: In \cite{HuyTalk}, Huybrechts explains that this result was mentioned by André Weil in \cite{Weil79} without attribution, but may originally be due  to Andreotti. 
    More recently, the problem has been explored in greater generality: For example,  \cite[Theorem 1']{LiWi77} establishes cohomological criteria for the local Torelli theorem to be satisfied and \cite[Proposition 3.14]{HuyGenCY} proves it for generalized Calabi-Yau manifolds. The results in both texts certainly apply to complex tori. 
    %$X=\Cnum^2/\Lambda$.
\end{rmk}

By Theorem \ref{thm:loctor},  locally near $I$, the moduli space of complex structures in these cases is a complex manifold contained in the projectivization  of the 2nd cohomology. In particular, changing the holomorphic symplectic form within its cohomology class results in an isomorphic complex structure. 
 It is well-known that the domain of the period map $\cP$ is non-Hausdorff.

For both 4-tori and K3 surfaces -- i.e. the manifolds considered in Theorem \ref{thm:loctor} -- there are also known results about the global properties of the period map. For the purpose of this text, we point out the following
(for the K3 surface see \cite[Thm. 0.4]{HuyTalk}\cite[\S 2.1]{Buchdahl}, for the 4-torus we refer to the beginning of \cite[\S 2.3]{Buchdahl}): 
\begin{thm}  
\label{thm:globtor}
Let $(M,I)$ be either the 4-torus or the K3-manifold, endowed with any complex structure inducing the standard orientation.
Then the period map \eqref{eq:globgtor} is   surjective. 
\end{thm}

\begin{rmk}
A version of the Local Torelli theorem holds for compact hyperk\"ahler manifolds in general.  (\cite{BeauvilleLocTorelliHK}, see also \cite[Theorem 4.1 and \S3.1] {HuybrechtsHKBourbaki}). In this paper we will not address the question of whether our main result  - Theorem \ref{prop:CMHR}, which  relies on the Local Torelli Theorem - can be extended to compact hyperk\"ahler manifolds of higher dimensions.  
\end{rmk}

\section{Moduli spaces of spacefilling branes  on symplectic $4$-manifolds}  
\label{sec:defdiff}

This section contains our main results. 
Consider a compact symplectic 4-manifold  admitting a spacefilling brane such that the  complex structure is K\"ahler. {Then this complex structure admits a holomorphic symplectic form, by Proposition \ref{prop:charbrane}.}
This restricts the underlying smooth manifold to the 4-torus and the K3 manifold, as we have established in Remark \ref{rem:holsymK4}. 
Our results are geometric statements about the moduli space of spacefilling brane structures: the moduli space is a non-Hausdorff smooth manifold, it is non-compact, of dimension $4$ for the 4-torus and  $20$ for the K3 manifold, and it carries a natural Lorentzian metric. These results are presented in \S \ref{subsec:smoothmod} and \S \ref{subsec:globalproperties}.  The previous subsections \S \ref{sec:branes4} and 
\S \ref{subsec:exactdef} are devoted to a characterization of  spacefilling branes on general symplectic 4-manifolds, and to the definition of moduli space of spacefilling branes.
\bigskip

We outline briefly how we obtain these results. 
We take the following point of view:  deforming a brane structure $F$ means
deforming the complex structure $I:=\omega^{-1}\circ F$ within the class of complex structures $I'$ for which
$\omega\circ I'$ is skew, by Proposition \ref{prop:charbrane} (This deformation of course corresponds to a deformation of the 2-form $F$).

We make use of the properties of compact oriented $4$-manifolds we recalled in \S \ref{sec:complex4man}, in particular the correpondence between complex structures and lines of complex 2-forms.
Complex structures as above then correspond\footnote{Notice that this is consistent with the fact that deforming the brane structure $F$ can be also be viewed as deforming the holomorphic symplectic form $\Omega=F+i\omega$ (changing the complex structure with respect to which the new form is holomorphic)  while keeping its imaginary part $\omega$ fixed.}  to complex 2-forms with imaginary part $\omega$.
Our statements about the moduli space of spacefilling branes are then obtained upon modding out by suitable equivalences and applying the Local Torelli theorem.

\subsection{Spacefilling branes on symplectic 4-manifolds}\label{sec:branes4}

In this subsection we restrict the correspondence recalled in \S \ref{subsec:4man} by  fixing the imaginary part of the complex 2-form.
In the next proposition we recover the characterization of spacefilling branes in symplectic 4-manifolds first observed by Kapustin-Orlov \cite[\S3]{KO}. 

\begin{prop}
\label{prop:OmegaF}
Let $M$ be an oriented 4-manifold, and fix a symplectic form $\omega$ inducing the given orientation. 
The conditions for a 2-form $F\in \Omega^2(M,\RR)$ to  define a spacefilling brane structure on $(M,\omega)$ are precisely:
 \begin{SB}
\item $F\wedge F= \omega \wedge \omega$ \label{list:pos}
\item $F\wedge \omega =0$ \label{list:orth}
\item $\dx F=0$. \label{list:cl}
\end{SB}
\end{prop}

\begin{proof}
 We know that $F\in \Omega^2(M,\RR)$ defines a spacefilling brane structure if{f}   the complex 2-form $\Omega=F+i\omega$ (with the fixed symplectic form $\omega$ as the imaginary part) is holomorphic for some complex structure on $M$ defining the given orientation, by Proposition \ref{prop:charbrane} and Remark \ref{rem:holDarboux}. This happens precisely when
$\Omega=F+i\omega$ satisfies the conditions HS(1)-HS(3) listed in Proposition \eqref{prop:HS}.

These   conditions   are equivalent to the above three equations for $F$. Indeed, the first two conditions are equivalent to $\Omega\wedge\Omega=0$, i.e. to $(HS 2)$. 
Using this, we have $\Omega\wedge \bar{\Omega}=2\omega\wedge \omega$, thus the condition    $(HS 1)$ is automatically satisfied.
Finally, since $\omega$ is closed, the third condition in the statement is equivalent to $(HS 3)$.
\end{proof} 

\begin{rmk}\label{rmk:or}
Notice that 
%since  $\omega$ is symplectic, 
condition SB(1) above implies that $F$ is non-degenerate and  that it induces the given orientation on $M$.
\end{rmk}

Now let $\omega$ be a symplectic form and $F$ a spacefilling brane structure, as in Proposition \ref{prop:OmegaF}. 
We can determine what its \emph{deformations} $F\mapsto F+\alpha$ with $\alpha\in\Omega^2(M,\RR)$ look like, by applying conditions \ref{list:pos}-\ref{list:cl} to $F+\alpha$.

\begin{lem}
    \label{lem:alphabrane}
 {Given a spacefilling brane $F\in\Omega^2(M,\RR)$, a} 2-form $F+\alpha$ defines a space-filling brane  on $(M,\omega)$ if{f} 

\begin{enumerate}[label=(\roman*)] 
\item $ F\wedge \alpha + \frac{1}{2}\alpha\wedge \alpha=0$ \label{list:defo1}
\item $\omega\wedge \alpha =0$ \label{list:defo2}
\item $\dx \alpha =0$
\end{enumerate}
\end{lem}

By linearising the equations  for $\alpha$ in Lemma \ref{lem:alphabrane}, we find the equations characterizing \emph{infinitesimal deformations} of the spacefilling brane $F$.
We need some preparation in order to express them in terms of the complex structure $I$ associated to the brane structure $F$.

 %We can re-express them thanks to the following lemma.  
\begin{lem}
\label{cor:11}
For any 2-form $\alpha\in \Omega^2(M,\RR)$, the following two statements are equivalent: 
\begin{itemize}
\item[a)] $F\wedge \alpha=0$ and $\omega\wedge \alpha =0$. 
\item[b)] $\alpha \in \Omega^{1,1}_{\RR}(M)$. 
\end{itemize} 
\end{lem} 
\begin{proof}

Notice that
 $F$ and $\omega$ are both elements of $(\Omega^{2,0}\oplus \Omega^{0,2})_{\RR}$  as a consequence of Proposition \ref{prop:charbrane}. Furthermore they are  everywhere non-degenerate, and they are orthogonal according to item \ref{list:orth}. 
 Since the  vector bundle $(T^*)^{(2,0)+(0,2)}_{\RR}$  has rank $2$, it follows that  $F,\omega$ form a  frame for it.  
By Remark \ref{rmk:ortho} we infer that any two-form $\alpha$ that is orthogonal to both $F$ and $\omega$ has to lie in $\Omega^{1,1}_{\RR}(M)$, and vice versa. 
\end{proof} 
\begin{cor}\label{cor:infdef}

The infinitesimal deformations of spacefilling branes are given by closed  {real} (1,1)-forms.
\end{cor}

Further, for  any closed $(1,1)$-form $\alpha$ that is isotropic with respect to the wedge product,   $t\alpha$  is an actual  deformation for all $t\in \RR$. This follows from Lemma \ref{lem:alphabrane}.

\subsection{The moduli space of spacefilling branes}
\label{subsec:exactdef}

Let $(M,\omega)$ be a compact symplectic manifold. In this section we introduce the notion of equivalence between spacefilling branes that we will use in the rest of this note.

Denote by $\operatorname{Symp}_0(M,\omega)$ the connected component of the symplectomorphism group that contains the identity. A possible notion of equivalence, which is natural from a geometric perspective, is to consider two  brane structures $F, F'\in \Omega^2(M,\RR)$  to be equivalent if they are related by a symplectomorphism $\phi \in \operatorname{Symp}_0(M)$, i.e. if 
    $\phi^{*}(F)=F'$. Notice that in that case we have $F'=F+\dx \xi$ for some $1$-form $\xi$, i.e. $F,F'$ lie in the same cohomology class. Conversely, Moser's theorem (applied to $F$) tells us that a sufficiently small deformation $F\mapsto F+\dx \xi$ to a new brane structure can always be achieved via a smooth isotopy, though we do  not know whether this preserves the symplectic form $\omega$.

In this text, we choose the potentially coarser equivalence relation induced by the closely related 
group of symplectomorphisms inducing the identity in cohomology:  
$$\Symp_*(M,\omega):=\Symp(M,\omega)\cap\Diff_*(M),$$
where $\Diff_*(M)$ was defined in \eqref{eq:diff*}. 
Our reason for {this choice of equivalence relation is based on the Local Torelli Theorem:}  by applying the Local Torelli Theorem \ref{thm:loctor} in the next subsections we will obtain results about the moduli space $\cM_{\omega}$, rather than  about the moduli space arising from $\Symp_0(M,\omega)$.

\begin{rmk}
While it is clear that 
 $\Diff_0(M)\subset \Diff_*(M)$ and 
$\Symp_0(M,\omega)\subset \Symp_*(M,\omega)$, it is not known whether  these inclusions are  strict. 
\end{rmk}

\begin{defn}\label{def:equibrane}
Two spacefilling brane structures $F,F'$ on $(M,\omega)$ are considered \emph{equivalent branes},  denoted $F\sim F'$, if 
\[ F'=\phi^* F \text{ for } \phi\in \Symp_*(M,\omega).\]
The \emph{moduli space of spacefilling brane structures on $(M,\omega)$} is 
\[ \mathcal{M}_{\omega} := \{F\in \Omega^2(M,\RR) \text{ spacefilling brane for }(M,\omega)\}/\sim. \]
\end{defn}

Notice that we have
\[ F\sim F' \Rightarrow 
%F'=F + \dx \xi, \xi \in \Omega^1(M),
[F]=[F'],
\]
but the converse is in general unknown (the Moser theorem applies only if $F$ and $F'$ are sufficiently close to each other).
The equivalence relation given by $F \sim_{H^2} F' :\Leftrightarrow [F]=[F'] \in H^2(M,\RR)$ is thus in general coarser than equivalence of branes as defined above. 

The remainder of this section focuses on learning about $\mathcal{M}_{\omega}$ via the cohomology $H^2(M,\RR)$, and establishing contexts in which the equivalence of branes $\sim$ and $\sim_{H^2}$ are the same. 
 
\subsection{Smoothness of the moduli space  of spacefilling branes on Kähler $4$-manifolds}
\label{subsec:smoothmod}

In this subsection $(M,\omega)$ denotes a compact symplectic 4-manifold admitting a spacefilling brane structure  whose corresponding complex structure is  Kähler. We fix the orientation of $M$ induced by $\omega$. 

The aim of this subsection is to  {prove} Theorem \ref{prop:CMHR}: The moduli space of branes is smooth. To do so, we first restrict the  period   map $\cP$ defined in \eqref{eq:globgtor}  in a suitable manner.

\subsubsection{Smoothness of the restricted period map}\label{subsec:restr}

As a consequence of   Proposition \ref{prop:charbrane}, the bijective correspondence of Proposition \ref{prop:corr1}  restricts to one between 
\begin{enumerate} [label=(\roman*)]
    \item  complex structures $I'$ such that $\omega\circ I'$ is skew\footnote{Notice that $I'$ defines the given orientation on $M$, by the text before Remark \ref{rmk:ortho} and by Proposition \ref{prop:OmegaF}.},
    
    \item  complex lines of 2-forms $\CC \Omega'$ in $\Omega^2(M,\Cnum)$   satisfying HS(1)-HS(3) and so that $Im(\Omega')=\omega$.
\end{enumerate}
 This correspondence reads 
 \begin{equation}\label{eq:ItoOmega}
  I'\mapsto \CC(\omega\circ  I'+i\omega).   
 \end{equation}

Therefore,  restricting the period
   map  $\cP$ in \eqref{eq:globgtor}
   %in \S \ref{subsec:localtorelli}
    one obtains a well-defined map 
\begin{equation}  \label{eq:globgtoromega}
\cP_{\omega}\colon 
\left\{
\begin{matrix}
\text{Complex structures $I'$ on $M$}\\
\text{such that $\omega\circ I'$ is skew}\\
%\text{defining the given orientation}
\end{matrix}
\right\}
/ \Diff_*(M)\to 
\{\Cnum[\Omega']\in \mathcal{Q} \text{ s.t. } Im[\Omega']=[\omega]\},
\end{equation}
which we will call \emph{restricted period map}.
Here $\mathcal{Q}\subset \mathbb{P}H^2(M,\Cnum)$ is the period domain, which we defined as the projectivization of $\QQQ\subset H^2(M,\Cnum)$ as in \eqref{eq:Q}.
The domain of the map is a quotient by which two complex structures are identified if they are related by an element of $\Diff_*(M)$. 

For later use, we state the following lemma.
\begin{lem}\label{lem:1}
Let $[\Omega']\in \QQQ$. Then the only complex multiple of $[\Omega']$ with imaginary part $Im([\Omega'])$ is $[\Omega']$ itself.
\end{lem}
\begin{proof}
Write $[\Omega']=[F']+i[\omega']$. The fact that $[\Omega']\in \QQQ$ can be written as 
\[ [\omega']\wedge [F'] =0,\quad [F']\wedge [F']=[\omega']\wedge[\omega'] >0.\]
This implies that $[F']$ and $[\omega']$ are linearly independent elements of the real vector space $H^2(M,\RR)$.

The imaginary part of $(a+ib)[\Omega']$ is 
%Let $a+ib$ be a complex number. 
$b[F']+a[\omega']$, hence it equals  $[\omega']$ only when $a+ib=1$.
\end{proof}
 
We restrict the one-to-one correspondence  
of the Local Torelli theorem to the classes of complex structures corresponding to brane structures (w.r.t. $\omega$). Doing so, it follows that the restricted period map $\cP_{\omega}$ of
\eqref{eq:globgtoromega} is locally a bijection (cf. the text just after Theorem \ref{thm:loctor}).

\begin{lem}\label{lem:3equiv}
Let $I$ be a complex structure on $M$ such that $\omega\circ I$ is skew, and so that $I$ is  K\"ahler.
Restricting the period map $\cP$ yields a bijection between:
\begin{itemize}
    \item[a)] small deformations $I'$ of the complex structure $I$ such that $\omega\circ I'$ is skew, up to diffeomorphisms in $\Diff_*(M)$ 
    \item[b)] small deformations $\Cnum[\Omega']$ of the complex line $\Cnum[\Omega]$ in $\mathcal{Q}$ such that $Im[\Omega']=[\omega]$.
\end{itemize}
\end{lem}
\begin{proof}
Applying the Local Torelli Theorem \ref{thm:loctor} one obtains an injective map from \emph{a)} to \emph{b)}. It actually provides a bijection \emph{a)-b)}. Indeed,    given a line as in \emph{b)}, we claim that its   unique element $A$  with imaginary part $[\omega]$ has a representative  satisfying HS(1)-HS(3) in Proposition \ref{prop:HS}  whose imaginary part is precisely $\omega$. The complex line spanned by this representative hence lies in the image of the correspondence recalled in \eqref{eq:ItoOmega}. 

To show the claimed existence of the representative of $A\in H^2(M,\Cnum)$, we argue as follows: 
There exists a smooth universal family of deformations $X\to S$ of $(M,I)$, parametrized by an open neighborhood $S$ of $\CC[\Omega]$ in the period domain $\mathcal{Q}$. For K3 surfaces this follows from  \cite[\S 6.2.3, Corollary 2.7 and Proposition 2.8]{HuyLecturesK3}, and for 4-tori from 
 \cite[Thm. VII.1, Thm. I. 52, Ex. I.54]{ManettiLectDefComplexMan}.
Denote by $I'$  the complex structure on the fiber over $\CC[\Omega']$;
then $I'$ lies nearby $I$, and $I'$ is a representative of the class of complex structures that corresponds to $\CC[\Omega']$
under Theorem \ref{thm:loctor}.
Let $\Omega'$ be the unique 2-form lying in the complex line of forms associated to $I'$ (as in Proposition \ref{prop:corr1}) and whose imaginary part in cohomology is $[\omega]$ (the uniqueness is ensured by Lemma \ref{lem:1}).
Then $\Omega'$  is a representative of  $A$  lying nearby $\Omega=F+i\omega$ and satisfying HS(1)-HS(3). In particular, writing $\Omega'=F'+i\omega'$, we have that $\omega'$ lies nearby $\omega$ and lies in the same cohomology class, so by Moser's theorem there is $\varphi\in \Diff_*(M)$ such that $\varphi^*\omega'=\omega$. The desired representative of $A$ is   $\varphi^*(\Omega')$.
\end{proof}

We now turn to the question of smoothness of $\cP_{\omega}$.
\begin{lem}\label{lem:trans}
The set $\{\Cnum[\Omega']\in \mathcal{Q} \text{ s.t. } Im[\Omega']=[\omega]\}$
is a submanifold
 of $\mathcal{Q}$.
\end{lem}
\begin{proof}
\textit{ Claim: $ \QQQ$ and $H^2(M,\Rnum)+i[\omega]$  are transverse submanifolds in $H^2(M,\Cnum)$.}

To see this, let $[\Omega']=[F']+i[\omega]$ lie in the intersection. Given an arbitrary vector $[F'']+i[\omega'']\in H^2(M,\Cnum)=T_{[\Omega']}H^2(M,\Cnum)$,
for any $a,b\in \RR$ we can write it as
$$\big([F''+aF'+b \omega]+i[\omega'']\big)-[aF'+b \omega].$$
The second summand is real, thus belongs to $T_{[\Omega']}(H^2(M,\Rnum)+i[\omega])$.
There exists a (unique) choice of $a,b$ such that the first summand lies in $T_{[\Omega']}\QQQ=[\Omega']^{\perp}$, as one computes using $[F']\wedge [\omega]=0$ and $[F']\wedge[F']>0$, $[\omega]\wedge[\omega]>0$ (which hold by the proof of Lemma \ref{lem:1}).

By the claim, the intersection $\QQQ \cap (H^2(M,\Rnum)+i[\omega])$ is a submanifold of $\QQQ$. It meets at most once every fiber of the submersion $p\colon \QQQ \to \mathbb{P}\QQQ =\mathcal{Q}$, as a consequence of Lemma \ref{lem:1}. Further, 
the tangent space at $[\Omega']$ of $\QQQ \cap (H^2(M,\Rnum)+i[\omega])$ intersects trivially the tangent space $\mathbb{C}[\Omega']$ to the $p$-fiber:   the only real element in the line $\mathbb{C}[\Omega']$ is the zero vector, since $[F']$ and $[\omega]$ are linearly independent.

Hence the image of $\QQQ \cap (H^2(M,\Rnum)+i[\omega])$ under the projection $p$,  which is the set in the statement of this lemma, is a submanifold of $\mathbb{P} \QQQ$.
\end{proof}

For the reader's convenience, we display again the period   map  $\cP$ from \eqref{eq:globgtor} and the restricted period map
$\cP_{\omega}$ from \eqref{eq:globgtoromega}.

 \begin{center}
\begin{tikzcd}
\cP\colon 
 \left\{
 \begin{matrix}
\text{Complex structures on $M$}\\
\text{defining the given orientation}
\end{matrix}
\right\}
/ \Diff_*(M)  \arrow[r ] & \mathcal{Q} \\

\cP_{\omega}\colon 
\left\{
\begin{matrix}
\text{Complex structures $I'$ on $M$}\\
\text{such that $\omega\circ I'$ is skew}\\
%\text{defining the given orientation}
\end{matrix}
\right\}
/ \Diff_*(M)
%\arrow[u, hookrightarrow]
\arrow[u, phantom, sloped, "\subset"]
\arrow[r] & 
\arrow[u, phantom, sloped, "\subset"]
\{\Cnum[\Omega']\in \mathcal{Q} \text{ s.t. } Im[\Omega']=[\omega]\} 
\end{tikzcd}  
\end{center}

\begin{cor}\label{cor:restrpermap}
Let $(M,\omega)$ be a compact symplectic 4-manifold admitting a spacefilling brane structure, and assume that the corresponding complex structure is K\"ahler.

The domain of the restricted period map $\cP_{\omega}$ is a submanifold of the domain 
%$\{\text{Complex structures on $M$ defining the given orientation}\}$ 
of the period  map  $\cP$.

Further, the restricted period map  $\cP_{\omega}$ is locally a diffeomorphism onto its image; in particular it is smooth, and thus continuous w.r.t. the subspace topologies.
\end{cor}
\begin{proof}
The period map $\cP$ in \eqref{eq:globgtoromega} is a local diffeomorphism, by the Local Torelli theorem  \ref{thm:loctor}. (The theorem applies because every complex structure on $M$ is K\"ahler, see Remark \ref{rem:holsymK4}.)

The   restricted period map  $\cP_{\omega}$ is locally a bijection by Lemma \ref{lem:3equiv}, and its codomain is a submanifold of $\mathcal{Q}$
by Lemma \ref{lem:trans}.  Since being a submanifold is a local property, this proves the first statement. The second statement is a direct consequence.
\end{proof}

\subsubsection{Smoothness of the moduli space of branes}\label{sec:Q}

We start reformulating the \emph{codomain} of the   of the restricted period map $\cP_{\omega}$ of \eqref{eq:globgtoromega}.

\begin{defn}\label{def:Qw}
  We define        $\mathcal{Q}_{[\omega]}$  to be
    the subset of elements $[F']$ in $H^2(M,\RR)$ cut out by the conditions 
    \begin{equation}\label{eq:QR}
            [F']\wedge [\omega]= 0 ,\quad [F'] \wedge [F']=[\omega] \wedge [\omega].  
    \end{equation}  
\end{defn}

 Notice that $\mathcal{Q}_{[\omega]}$ is a real quadric  submanifold  of $H^2(M,\Rnum)$ of codimension $2$. (This can be seen using the regular value theorem\footnote{Indeed, $\mathcal{Q}_{[\omega]}$ is the preimage under the map
$$\phi\colon [\omega]^{\perp}\to H^4(M,\Rnum)\cong \RR, \quad [F']\mapsto [F']\wedge [F']$$ of the value $[\omega] \wedge [\omega]$. This is a regular value since for every $[F']\in \mathcal{Q}_{[\omega]}$ we have 
 $(d_{[F']}\phi)([F'])=\left.\frac{d}{dt}\right\vert_{0}\phi([F']+t[F'])=2[F']\wedge [F']=2[\omega] \wedge [\omega] \neq 0$.}, but also from Proposition \ref{thm:Qw} later on.)

\begin{lem}\label{rem:qomega}
There is a diffeomorphism
$$\{\Cnum[\Omega']\in \mathcal{Q} \text{ s.t. } Im[\Omega']=[\omega]\}\;\cong\; \mathcal{Q}_{[\omega]},$$
mapping the line $\Cnum[\Omega']$ to $Re([\Omega'])$, where $[\Omega']$ is determined by $Im[\Omega']=[\omega]$. 
\end{lem} 
\begin{proof}
This map is a diffeomorphism since it is obtained composing the diffeomorphisms

\[ \{\Cnum[\Omega']\in \mathcal{Q} \text{ s.t. } Im[\Omega']=[\omega]\} \longleftarrow  \QQQ \cap (H^2(M,\Rnum)+i[\omega])\longrightarrow \mathcal{Q}_{[\omega]}, \]
where the left map is the restriction of the projection $p$ as in the proof of Lemma \ref{lem:trans}. The right map is the restriction of the isomorphism $(H^2(M,\Rnum)+i[\omega])\cong H^2(M,\Rnum)$ which simply takes the real part. The right map is a bijection by the proof of Lemma \ref{lem:1}.
\end{proof}

Notice that upon the identification of Lemma  \ref{rem:qomega},  the restricted period map $\cP_{\omega}$ of 
\eqref{eq:globgtoromega} reads
\begin{equation}\label{eq:restrictedcodomain}
    I' \text{ mod } \Diff_*(M,\omega) \mapsto [\omega\circ I'],
\end{equation}
as a consequence of \eqref{eq:ItoOmega}.

We now reformulate the \emph{domain}    of the restricted period map, by showing that it is   unchanged if we replace \enquote{$\Diff_*(M)$} with \enquote{$\Symp_*(M,\omega)$}.

\begin{lem} \label{prop:equiv}
Let $(M,\omega)$ be a compact symplectic 4-manifold. 
Let $F'$ be a brane structure, with corresponding complex structure  $I':=\omega^{-1}\circ F'$, and $F''$ a second nearby brane structure with complex structure $I''$. 

If $I'$ and $I''$ are related by a diffeomorphism $\phi$ lying in $\Diff_*(M)$,
then $\phi$ is a symplectomorphism with respect to $\omega$. In other words, $\phi$ lies in $\Symp_*(M,\omega)$.
\end{lem}
 
Note that for the above proposition there is no requirement on  $M$ being Kähler.

\begin{proof} 
Recall that by Proposition \ref{prop:charbrane} $F'+i\omega$ is holomorphic symplectic w.r.t. $I'$, and similarly $F''+i\omega$ is holomorphic symplectic w.r.t. $I''$.
Since $I'$ is related to $I''$ by the diffeomorphism $\phi$, it follows that $\phi^*(F''+i\omega)$ is holomorphic symplectic with respect to $I'$. By Proposition \ref{prop:corr1}, $F'+i\omega$ and $\phi^*(F''+i\omega)$ must be holomorphic sections of the same line bundle of 2-forms, i.e. 
\[\phi^*(F''+i\omega)=\lambda(F'+i\omega) \]
for some $\lambda \in \Cnum\setminus\{0\}$.
Using that $\phi^*$ is the identity on cohomology, we obtain
\[[F'']+i[\omega]=\lambda([F']+i[\omega]). \] Lemma \ref{lem:1} implies that $\lambda=1$ and thus $\phi^*\omega =\omega$.
\end{proof}

Recall from \S \ref{subsec:exactdef} that we denote by $\cM_{\omega}$ the moduli space of spacefilling brane structures on $(M,\omega)$ up to diffeomorphisms in $\Symp_*(M,\omega)$.
%, as introduced in Def. \ref{def:equibrane}.
Applying to the  domain   of the restricted period map the correspondence $I'  \mapsto \omega\circ I'$ of Proposition \ref{prop:charbrane}, and using Lemma \ref{prop:equiv}, we obtain an identification between the domain of the restricted period map $\cP_{\omega}$ of \eqref{eq:globgtoromega}  and $\cM_{\omega}$. In conclusion, we obtain an identification of  the restricted period map $\cP_{\omega}$ of \eqref{eq:globgtoromega} with the map
 \begin{align}\label{eq:globtor}
\Phi \colon \cM_{\omega}&\to \mathcal{Q}_{[\omega]} \\
 F' \text{ mod } \Symp_*(M,\omega) &\mapsto [F'], \nonumber
\end{align}
as can be seen using \eqref{eq:restrictedcodomain}.

 Applying Corollary \ref{cor:restrpermap} we obtain:

\begin{thm}\label{prop:CMHR}
Let $(M,\omega)$ be a compact symplectic 4-manifold admitting a spacefilling brane structure, and assume that the corresponding complex structure is K\"ahler.

Then the moduli space of brane structures $\cM_{\omega}$  is a non-Hausdorff smooth manifold. 

Further, the map $\Phi$ of \eqref{eq:globtor} is smooth, and is a local diffeomorphism.
\end{thm}
{Recall  that the smooth manifold  $M$  in this theorem must be either the K3 manifold or the 4-torus, by  Proposition \ref{prop:charbrane} and Remark \ref{rem:holsymK4}.}

In particular, for any spacefilling brane structure $F'$, 
the map $\Phi$
yields a  diffeomorphism between:

\begin{enumerate} [label=(\alph*)]
\item  a neighborhood of the class of the spacefilling brane $F'$ in  $\cM_{\omega}$,

\item small deformations of $[F']$ in  $\mathcal{Q}_{[\omega]} \subset H^2(M,\RR)$.
\end{enumerate}

\subsubsection{Remarks on the map $\Phi$}\label{sec:remarksPhi}

For any spacefilling brane structure $F$, the translation by $-[F]$ maps $\mathcal{Q}_{[\omega]}$ to
$\mathcal{Q}_{[\omega]}^{[F]}$, the subset of  $[\alpha]$'s in  $H^2(M,\Rnum)$ cut out by the following equations: 
\begin{align}
[\omega]\wedge [\alpha] &=0 \label{eq:defo1} \\
[F]\wedge [\alpha] +\frac{1}{2}[\alpha]\wedge[\alpha] &=0 \label{eq:defo2}.
\end{align}
Theorem \ref{prop:CMHR} thus states that the moduli space $\cM_{\omega}$ around  $[F]$ is diffeomorphic to  a neighborhood of the origin in $\mathcal{Q}_{[\omega]}^{[F]}$, via the assignment $(F+\alpha \text{ mod } \Symp_*(M,\omega))\mapsto [\alpha]$.

\begin{rmk}\label{rmk:refine}
\textbf{(Refining Corollary \ref{cor:infdef})}
Since the inner product induced on $H^2(M,\Rnum)$ by the wedge product on a compact Kähler manifold is in fact a \emph{Hodge structure}, we can now make a   statement about infinitesimal deformations, stronger than the one of Corollary \ref{cor:infdef},  based on the  local  Torelli theorem: 

The linearised versions of equations (\ref{eq:defo1}), (\ref{eq:defo2}) are 
\begin{align*}
[\omega]\wedge [\alpha] &=0 \\ 
[F] \wedge [\alpha] &=0.
\end{align*}
This stipulates precisely that an infinitesimal deformation modulo trivial deformations by symplectomorphisms, $[\alpha]$, has to lie in the orthogonal subspace to the space $\left(H^{2,0}(M)\oplus H^{0,2}(M)\right)_{\Rnum}$ defined in Appendix \ref{sec:dRD}.
The latter is $H^{1,1}(M)_{\Rnum}$
by the orthogonal decomposition \eqref{eq:Hodge} into types w.r.t. the complex structure $I:=\omega^{-1}\circ F$,
and by Corollary \ref{cor:realsplitting}.
In other words, the tangent space to $\cM_{\omega}$ at $[F]$ is $H^{1,1}(M)_{\Rnum}$, refining Corollary \ref{cor:infdef}.
\end{rmk}

We continue denoting by $(M,\omega)$  a compact symplectic 4-manifold admitting a spacefilling brane structure whose corresponding complex structure is  Kähler.
In general, although the period map \eqref{eq:globgtor} is surjective (see Theorem \ref{thm:globtor}),  we do not know whether the map $\Phi$ of \eqref{eq:globtor} is surjective. In the case of the 4-torus however we can check the surjectivity explicitly.

\begin{prop}\label{prop:torussurj}
 Let $\omega$ be any symplectic form on the 4-torus  $T^4=\RR^4/(2\pi\ZZ)^4$  admitting a spacefilling brane structure. Then the map  
  \begin{align*}
        \Phi \colon \cM_{\omega}&\to \mathcal{Q}_{[\omega]} \\
        F' \text{ mod } \Symp_*(M,\omega) &\mapsto [F']. 
        \end{align*}
of \eqref{eq:globtor} is surjective.
\end{prop}

\begin{proof}
We start establishing the following

\textit{Claim: the symplectic form $\omega$ is diffeomorphic to a constant   symplectic form on $T^4$.}

To see this, recall that   any complex structure on $T^4$ is related by a diffeomorphism to the standard complex structure on  $\RR^4/\Lambda$ (i.e. the one induced 
from the  standard complex structure on $\RR^4=\CC^2$), for some full rank lattice $\Lambda$ in $\RR^4$.
In turn, there is a linear automorphism of $\RR^4$ which maps $(2\pi\ZZ)^4$ to $\Lambda$, and hence induces a diffeomorphism $T^4\to \RR^4/\Lambda$. Under this diffeomorphism, 
the standard complex structure on  $\RR^4/\Lambda$ corresponds to a constant complex structure $I'$ on $T^4$. Further,  the  1-forms on  $\RR^4/\Lambda$ induced by the coordinate 1-forms on $\RR^4$ correspond to constant 1-forms on $T^4$, and therefore 
the holomorphic symplectic form 
%for the standard complex structure 
on $\RR^4/\Lambda$ induced from $dz_1\wedge dz_2$ on $\CC^2$ corresponds to a constant form $\Omega'$  on $T^4$, which is holomorphic symplectic for $I'$. 

Denote by $F$ a spacefilling brane structure for $\omega$; it gives rise to a complex structure $I:=\omega^{-1}\circ F$. By the above, there is a diffeomorphism $\phi$ of $T^4$ mapping $I$ to  a constant complex structure $I'$, which furthermore admits a constant holomorphic symplectic form $\Omega'$.
By the uniqueness up to scalar multiples in Proposition \ref{prop:corr1}, and by Proposition \ref{prop:charbrane}, this means that $(\phi^{-1})^*(F+i\omega)$ is a constant multiple of  $\Omega'$. In particular $(\phi^{-1})^*(\omega)$ is a  constant 2-form, proving the claim.

To finish the proof, it suffices to show:

\textit{Claim:  When $\omega$ is a \emph{constant}   symplectic form on $T^4$, the map \eqref{eq:globtor} is surjective.}

 Consider the 1-forms $$\xi_1,\eta_1,\xi_2,\eta_2\in \Omega^1(T^4)$$   whose pullback to $\CC^2=\RR^4$ equal $\dx x_1,\dx y_1,\dx x_2,\dx y_2$ respectively. These 1-forms are constant, and their cohomology classes provide a basis of $H^1(T^4,\RR)$. Since the latter, together with $H^0(T^4,\RR)$, generates the whole cohomology, 
it follows 
that taking the cohomology class provides an \emph{isomorphism of graded algebras} from the   constant forms on $T^4$ to $H^{\bullet}(T^4,\RR)$.
In particular, a constant 2-form  $F'$ satisfies the first two equations in  Proposition \ref{prop:OmegaF} (i.e. it is a brane for $\omega$) if{f} its cohomology class $[F']$ satisfies \eqref{eq:QR} (i.e. it lies in $\mathcal{Q}_{[\omega]}$). 
\end{proof}

 The following general remark relates  classes  in $\Qw$ with the existence of representatives which are spacefilling branes.
 \begin{rmk}\label{rem:classbrane}
Given a symplectic form $\omega$, consider a class $A\in \Qw$.

\begin{itemize}
    \item[i)] 
In general, all we know is that there exists \emph{some}  symplectic form $\omega'$ cohomologoues to $\omega$ and admitting a spacefilling brane lying in $A$.

Indeed,
under the diffeomorphism of  Lemma \ref{rem:qomega}, the element $A\in \Qw$ corresponds to an element of  $\mathcal{Q}$. Since the period map \eqref{eq:globgtor} is surjective by Theorem \ref{thm:globtor}, there exists a complex structure $I$ compatible with the given orientation admitting a holomorphic symplectic form $F'+i\omega'$ with $[F']=A$ and $[\omega'] =[\omega]$. The symplectic form $\omega'$ therefore admits a brane whose class is $A$, namely $F'$ (see Proposition \ref{prop:charbrane}). 

  \item[ii)] If $A$ belongs to the image of the  map in \eqref{eq:globtor}
\begin{align*}
\Phi \colon \cM_{\omega}&\to \mathcal{Q}_{[\omega]} \\
F' \text{ mod } \Symp_*(M,\omega) &\mapsto [F'], 
\end{align*}
then by construction there is at least one representative $F$ of $A$ that is a spacefilling brane for $\omega$. The converse also holds.
\end{itemize}     
\end{rmk}

\subsection{Global properties of the moduli space of branes  on Kähler $4$-manifolds}
\label{subsec:globalproperties}

For this subsection,  as in the previous one,  let $(M,\omega)$ be a compact symplectic 4-manifold admitting a spacefilling brane structure, such that the corresponding complex structure is K\"ahler. As we saw in Remark \ref{rem:holsymK4}, there are just two classes of such objects: complex tori and K3 surfaces.

 We want to display geometric properties of the moduli space of branes $\cM_{\omega}$. We first focus on the geometry of the smooth quadric  $\Qw\subset H^2(M,\Rnum)$ of Definition \ref{def:Qw}, since $\cM_{\omega}$ is closely related to the latter
 by Theorem \ref{prop:CMHR}.

\begin{prop} \label{thm:Qw}
    The quadric $\Qw\subset H^2(M,\Rnum)$ is diffeomorphic to a cylinder $S^1\times \Rnum^{b_2-3}$, where $b_2=\dim(H^2(M,\Rnum))$. 

    Furthermore, the inner product on $H^2(M,\Rnum)$ restricts on the tangent spaces of $\Qw$ to a 
   Lorentzian metric of signature $(1,b_2-3)$. 
\end{prop}

In particular the topology and smooth structure of $\Qw$ are independent of the class of the symplectic form $[\omega]$.

\begin{proof}
    Recall that for any K\"ahler complex structure $I$, 
    the complex  de Rham cohomology decomposes as
    %there is a  decomposition into Dolbeault cohomology groups
    \[ H^2(M,\Cnum)= H^{2,0}(M) \oplus H^{1,1}(M) \oplus H^{0,2}(M).\]
    By Corollary \ref{cor:realsplitting} this induces an orthogonal decomposition of real de Rham cohomology
    \[ H^2(M,\Rnum) = (H^{2,0}\oplus H^{0,2})_{\Rnum}(M) \oplus H^{1,1}_{\Rnum}(M)\]
    with respect to the bilinear form $\wedge$. 
Here $(H^{2,0}\oplus H^{0,2})_{\Rnum}(M)$ coincides with the $\dx$-closed  elements of $(\Omega^{2,0}\oplus \Omega^{0,2})_{\Rnum}$ modulo its exact elements, and similarly for $H^{1,1}_{\Rnum}(M)$. (See Appendix \ref{sec:dRD}.)
    
Take an arbitrary spacefilling brane $F$. Denote by $I$  its associated complex structure.
    As established in \S\ref{sec:branes4}, we have $[\omega],[F]\in \HtR$ and $[\omega]\wedge[\omega]=[F]\wedge [F] > 0$. In fact, $[\omega]$ and $[F]$ form an orthogonal basis for $\HtR$, which is thus a positive-definite subspace of $H^2(M,\Rnum)$. 
    Its orthogonal complement $H^{1,1}_{\Rnum}(M)$ has signature $(1,b_2-3)$, since the signature of the 
inner product on $H^2(M,\RR)$ is $(3,b_2-3)$.
    Pick an element $b\in H^{1,1}_{\Rnum}(M)$ with $b\wedge b= [\omega]\wedge [\omega]$ and denote its orthogonal complement $\langle b\rangle^{\perp}\subset H^{1,1}_{\Rnum}(M)$ by $H^{1,1}_{<0}$. 

    Now, by equation (\ref{eq:QR}) any element $[F']$ of $\Qw$ has to be of the form 
    \[ [F'] = \lambda [F] + \eta b+ c,\]
    where $\lambda,\eta\in \Rnum$ and $c\in H^{1,1}_{<0}$, and satisfy 

      \[ (\lambda^2+\eta^2)([\omega]\wedge [\omega]) + c\wedge c= [\omega]\wedge [\omega].  \]

By choosing an orthogonal basis for $H^{1,1}_{<0}$ consisting of vectors squaring to $-[\omega]\wedge [\omega]$, we obtain an orthogonal basis with coordinates $y_4,\dots,y_{b_2}$ for $H^{1,1}_{<0}$. In these coordinates, the real quadric $\Qw$ reads as
    \[ \lambda^2 + \eta^2 - y_4^2 -\dots - y_{b_2}^2 = 1. \]
    The solution to this quadric is a smooth cylinder  diffeomorphic to $S^1\times \Rnum^{b_2-3}$. 

    The tangent space to $\Qw$ at $[F]$    is $H^{1,1}_{\Rnum}(M)$,  by Remark \ref{rmk:refine}, and we saw above that the signature there is $(1,b_2-3)$.
\end{proof}

\begin{rmk}
    While the smooth structure of $\Qw$ is independent of $[\omega]$, how this cylinder lies inside $H^2(M,\Rnum)$ obviously changes. As a consequence, the intersection $\Qw\cap 2 \pi H^2(M,\Znum)$ does depend on the class of the symplectic form. This is relevant when prequantised branes are considered, c.f. the companion paper \cite{PrequantisedSpacefillingBranes}.
    %{c.f. \S \ref{subsec:integral}.} 
\end{rmk}

In fact, we can make a stronger statement in specifying an isometry of $\Qw$ and $S^1\times \Rnum^{b_2-3}$ equipped with a particular Lorentzian metric:   

\begin{prop}\label{prop:metric}
The real quadric $\Qw$ equipped with the Lorentzian metric inherited from $H^2(M,\Rnum)$ is isometric to $S^1\times \Rnum^{b_2-3}$ equipped with the following metric in block form: 
\begin{align*}
    g = \left(\begin{matrix}
       \rho & 0 \\ 
       0 & \Gamma
    \end{matrix}\right)
\end{align*}
which is such that the metric $\rho$ on each copy of $S^1$ is positive, while the metric $\Gamma$ on each copy of $\Rnum^{b_2-3}$ is negative definite. 
\end{prop} 

The formulae for the metric $g$ provided in the proof make clear that  different classes $[\omega]$ give rise to metrics on $S^1\times \Rnum^{b_2-3}$  which are constant multiples of each other.

\begin{proof}
    We will specify a diffeomorphism 
    \[ \phi: S^1\times \Rnum^{b_2-3} \rightarrow \Qw \]
    and show that the preimage of the induced metric on $\Qw$ is as above. 

    Take coordinates $\lambda, \eta, y_4,\dots,y_{b_2}$ for the orthogonal complement of $[\omega]$ in $H^2(M,\Rnum)$ as in the proof of Proposition \ref{thm:Qw}, so that $\Qw$ is cut out by the equation 
    \[ \lambda^2 + \eta^2 - y_4^2 -\dots - y_{b_2}^2 = 1. \] 
  
    Take coordinates $(\theta,y_4,\dots,y_{b_2})$ for $S^1\times \Rnum^{b_2-3}$, with $\theta$ the angle on the circle and -- abusing  notation -- $y_4,\dots,y_{b_2}$ the standard coordinates on $\Rnum^{b_2-3}$.
      Denote: 
    \begin{align*}
    \bar{y}&:=(y_4,\dots,y_{b_2}) \in \Rnum^{b_2-3} \\
    r^2&:= \vert\vert\bar{y}\vert\vert^2=\sum_{i=4}^{b_2} y_i^2.
    \end{align*}
    
   Set  
    \begin{align*}
        \phi: S^1\times \Rnum^{b_2-3} &\rightarrow \Qw \\ 
        \phi(\theta,\bar{y}) &:= (\sqrt{1+r^2}\cos(\theta),\sqrt{1+r^2}\sin(\theta), \bar{y}).
    \end{align*}
    Note that $\phi$ is a smooth bijection onto its image, which is $\Qw$. With the following calculation and the implicit function theorem, we see that it is a diffeomorphism:  

    At a point $(\theta_0,\bar{y}_0)\in S^1\times \Rnum^{b_2-3}$, the tangent space is spanned by: 
    \begin{align*}
    \party{}{\theta} &=\left.\frac{\dx}{\dx t}\right\vert_{t=0}\gamma_{\theta}(t) = \left.\frac{\dx}{\dx t}\right\vert_{t=0}(\theta_0+t,\bar{y}_0) \\ 
    \party{}{y_i} &=\left.\frac{\dx}{\dx t}\right\vert_{t=0}\gamma_{y_i}(t)= \left.\frac{\dx}{\dx t}\right\vert_{t=0} \left(\theta_0,(y_4)_0,\dots,(y_i)_0+t,\dots,(y_{b_2})_0\right).
    \end{align*}
    Under $\phi$, these tangent vectors map to: 
    \begin{align*}
    \phi_*\left(\party{}{\theta}\right) &=\left.\frac{\dx}{\dx t}\right\vert_{t=0}(\phi\circ\gamma_{\theta})(t) 
    = \left.\frac{\dx}{\dx t}\right\vert_{t=0} \left(\sqrt{1+r_0^2}\cos(\theta_0+t),\sqrt{1+r_0^2}\sin(\theta_0+t), \bar{y}\right) \\ 
    &= -\sqrt{1+r_0^2}\sin(\theta_0)\party{}{\lambda} + \sqrt{1+r_0^2}\cos(\theta_0) \party{}{\eta} \\
    \phi_*\left(\party{}{y_i}\right) &=\left.\frac{\dx}{\dx t}\right\vert_{t=0}(\phi\circ\gamma_{y_i})(t) \\
    &= \left.\frac{\dx}{\dx t}\right\vert_{t=0} 
    \left(\sqrt{1+r_0^2+2t(y_i)_0+t^2}\cos(\theta_0),\sqrt{1+r_0^2+2t(y_i)_0+t^2}\sin(\theta_0),\dots,(y_i)_0+t,\dots\right) \\ 
    &= \frac{(y_i)_0}{\sqrt{1+r_0^2}} \cos(\theta_0)\party{}{\lambda} + \frac{(y_i)_0}{\sqrt{1+r_0^2}} \sin(\theta_0)\party{}{\eta} + \party{}{y_i}.
    \end{align*} 

    By construction, $\langle \party{}{\lambda},\party{}{\lambda}\rangle = \langle \party{}{\eta},\party{}{\eta}\rangle = [\omega]\wedge [\omega]$ and $\langle \party{}{y_i},\party{}{y_i}\rangle=-[\omega]\wedge [\omega]$, so we find that at each point $(\theta,\bar{y})\in S^1\times \Rnum^{b_2-3}$: 
    \begin{align*}
    \rho(\theta,\bar{y}):= g_{\theta,\theta}(\theta,\bar{y})&= \left\langle \phi_*\left(\party{}{\theta}\right),\phi_*\left(\party{}{\theta}\right)\right\rangle 
    =\left( \sqrt{1+r^2}\right) [\omega]\wedge[\omega] >0  \\ 
     g_{\theta,y_i}(\theta,\bar{y}) &= \left\langle \phi_*\left(\party{}{\theta}\right), \phi_*\left(\party{}{y_i}\right)\right\rangle 
    = 0 \\
    \Gamma_{y_i,y_i}(\theta,\bar{y}):= g_{y_i,y_i}(\theta,\bar{y}) &= \left\langle  \phi_*\left(\party{}{y_i}\right), \phi_*\left(\party{}{y_i}\right) \right\rangle 
    = \left(\frac{(y_i)^2}{1+r^2} - 1\right)[\omega]\wedge[\omega] <0 \\
    \Gamma_{y_i,y_j}(\theta,\bar{y}):= g_{y_i,y_j}(\theta,\bar{y}) &= \left\langle  \phi_*\left(\party{}{y_i}\right), \phi_*\left(\party{}{y_j}\right) \right\rangle 
    = \frac{y_i y_j}{1+r^2}[\omega]\wedge[\omega] \neq 0 \text{ for }i\neq j. 
    \end{align*}
    Since the  Lorentzian metric on  $\Qw$ has signature $(1,b_2-3)$ by Proposition \ref{thm:Qw}, the first two equations show that the metric does indeed split along the factorisation $S^1\times \Rnum^{b_2-3}$, with the factors being positive and negative definite, respectively.  The last two equations make the {negative-definite} metric $\Gamma$ explicit. 
\end{proof}

\begin{cor} \label{lem:compact}    For $M$ a K3  manifold  (respectively the 4-torus), endowed with any symplectic form $\omega$ admitting a spacefilling brane, the moduli space of brane structures $\cM_{\omega}$  
is not compact, and of  {real} dimension $20$ (respectively $4$). 
\end{cor}

\begin{proof}
The  {dimensions follow} from  {Theorem} \ref{prop:CMHR},  {Proposition} \ref{thm:Qw} and Remark \ref{rem:holsymK4}.

For the non-compactness, notice that the image of $\cM_{\omega}$ under the map \eqref{eq:globtor} is not compact. Indeed, if it was compact, it would be closed in   $\mathcal{Q}_{[\omega]}$ (since the latter is Hausdorff), but the image is also open by  {Theorem} \ref{prop:CMHR}. The image would then be the whole of the connected space $\mathcal{Q}_{[\omega]}$; the latter is not compact by  {Proposition} \ref{thm:Qw}, leading to a contradiction.  Now, if $\cM_{\omega}$ was compact, the  its image would also be.
\end{proof}

\subsubsection{An example: complex tori}

We now  take a detailed look at some explicit compact holomorphic symplectic Kähler manifolds, giving an explicit expression for the space $\Qw$ in terms of a fixed basis in cohomology. This provides an explicit example of the construction of Proposition \ref{thm:Qw}.

\begin{ex}\textbf{(Complex tori)} \label{ex:torus}
Set $M= T^4$. Any complex structure on $M$ is then obtained from a quotient $\Cnum^2/\Lambda$, where $\Lambda \subset \Cnum^2$ is a lattice of maximal rank, and  $\Cnum^2$ carries the standard complex structure.

Denote by 
\[\xi_1=\dx x_1,\ \eta_1=dy_1,\ \xi_2= \dx x_2,\ \eta_2=\dx y_2 \]
the canonical basis of one-forms on $\Cnum^2=\{(z_1=x_1+iy_1, z_2=x_2+iy_2)\}$. Under the fixed identification $T^4\cong \Cnum^2/\Lambda$, these descend to a basis of $H^1(M,\RR)$. The 2-form classes

\begin{equation}\label{eq:basisT4}
  [\xi_1]\wedge [\eta_1],\ [\xi_2]\wedge[\eta_2],\ [\xi_1]\wedge[\xi_2],\ -[\eta_1]\wedge[\eta_2],\ [\xi_1]\wedge [\eta_2],\ [\eta_1]\wedge [\xi_2]
\end{equation}
form a basis for $H^2(M,\RR)$. 

The holomorphic symplectic form on $M\cong \Cnum^2/\Lambda$ is descended from $\dx z_1 \wedge \dx z_2 $ on $\Cnum^2$, so it is given by 
\begin{align*}
\Omega&= (\xi_1+i\eta_1)\wedge(\xi_2+ i\eta_2) \\ 
&= \underbrace{\xi_1\wedge \xi_2 - \eta_1\wedge \eta_2}_{=:F} 
+i\underbrace{(\xi_1\wedge\eta_2 + \eta_1\wedge \xi_2)}_{=:\omega}.
\end{align*}

If we express a deformation in cohomology $[F]\mapsto [F]+[\alpha]$, where $[\alpha] \in H^2(M,\RR)$, in the basis \eqref{eq:basisT4} as a linear combination with real coefficients
\[ 
[\alpha]= f_1 [\xi_1]\wedge [\eta_1] + f_2 [\xi_2]\wedge [\eta_2] 
+ g_1 [\xi_1] \wedge [\xi_2] {-} g_2 [\eta_1] \wedge [\eta_2] 
+ h_1 [\xi_1] \wedge [\eta_2] + h_2 [\eta_1] \wedge [\xi_2], 
\]
 we find that the deformation equations \eqref{eq:defo1}, \eqref{eq:defo2} read as follows:
\begin{align*}
h_2&=-h_1 \\ 
{-} 2(g_1 g_2 {+} f_1 f_2 {+} h_1 h_2) + (g_1{+} g_2) 
&=0. \label{eq:tquadric}
\end{align*}
These equations are equivalent to the equation in 5 variables 
\begin{align*} {-}2(g_1 g_2 {+} f_1 f_2 {-}h_1^2) &+ (g_1{+} g_2) =0 \\ 
{g_1g_2 + f_1 f_2 -h_1^2} &{- \frac{1}{2}(g_1+g_2)=0.} 
\end{align*}

By changing coordinates, we find that this quadric is affinely equivalent to 
\[x_1^2{+}x_2^2 - x_3^2{-}x_4^2{-}x_5^2  =1, \] 
recovering the quadric equation for $\Qw$ from the proof of Proposition \ref{thm:Qw}.
A suitable coordinate transformation is: 
\begin{align*}
&x_1 = g_1+g_2 -1  
&x_2 = f_1+f_2 
&\phantom{mmmm} x_3 =2 h_1 & \\ 
&x_4 = g_1-g_2 
&x_5 = f_1-f_2 
&\phantom{x_5 = (f_1+f_2)}
\end{align*} 
\end{ex}

\appendix
\section{Real de Rham vs Dolbeault cohomology and induced orthogonal splitting}\label{sec:dRD}

Let $M$ be a compact K\"ahler manifold.
Recall from Remark \ref{rmk:ortho} that at the level of differential forms, it turns out to be useful to also split real forms by complex bidegree, for example for 2-forms: $\Omega^2(M,\Rnum)=(\Omega^{2,0}\oplus\Omega^{0,2})_{\Rnum}\oplus\Omega^{1,1}_{\Rnum}$. This splitting is orthogonal with respect to the pointwise inner product induced by $\wedge$. 

In Remark \ref{rmk:refine} and the
proof of Proposition \ref{thm:Qw}, 
we  consider a splitting of \emph{real cohomology} by complex bidegree, which for degree $(p,q)\oplus(q,p)$ can be defined in two ways: 
\begin{enumerate}
    \item 
Real elements of Dolbeault cohomology (i.e. those invariant under complex conjugation):
    
 $$\left(H^{p,q}\oplus H^{p,q}\right)_{\Rnum}:= \left\{a^{p,q}+\overline{a^{q,p}}\;\middle\vert\; a^{p,q} \in H^{p,q}_{\delb}(M)\right\}$$ 
    \item Quotients of real forms $$\tilde{H}^{(p,q)+(q,p)}_{\Rnum} := \frac{\left\{\dx\text{-closed  elements of } \left(\Omega^{p,q}\oplus \Omega^{q,p}\right)_{\Rnum}\right\}}{\left\{\dx\text{-exact elements of }\left(\Omega^{p,q}\oplus \Omega^{q,p}\right)_{\Rnum}\right\}}$$

\end{enumerate}
These two notions of cohomology turn out to be the same, as we prove in Lemma \ref{lem:HtildeH} for the relevant case $p+q=2$. The proof uses the $\del\delb$-Lemma in the following form, found in \cite{Huy}: 
\begin{lem}\textbf{(Corollary 3.2.10 in \cite{Huy})}\label{lem:HtildeH}
     Let $\alpha \in \Omega^{p,q}(M,\CC)$ be a $\dx\,$-closed form of type $(p,q)$ on a compact Kähler manifold. Then the following are equivalent: 
     \begin{enumerate}[label=(\roman*)]
         \item The form $\alpha$ is $\dx\,$-exact, i.e. $\alpha=\dx \beta$ for some $\beta \in  \Omega^{p+q-1}(M,\CC)$.
        
         \item The form $\alpha$ is $\del$-exact, i.e. $\alpha=\del \beta$ for some $\beta \in \Omega^{p-1,q}(M,\CC)$.
         \item The form $\alpha$ is $\delb$-exact, i.e. $\alpha=\delb \beta$ for some $\beta \in \Omega^{p,q-1}(M,\CC)$.
      \item The form $\alpha$ is $\del\delb$-exact, i.e. $\alpha=\del\delb\beta$ for some $\beta \in \Omega^{p-1,q-1}(M,\CC)$.
     \end{enumerate}
\end{lem}

 We can now prove:
 
\begin{lem}\label{lem:HRdclosed}
    For $p+q=2$, 
    \[ \left(H^{p,q}\oplus H^{q,p}\right)_{\Rnum} = 
    \tilde{H}^{(p,q)+(q,p)}_{\Rnum}
    \]
\end{lem}

\begin{proof}
We consider two cases.
\begin{enumerate} 
\item 
    Take a $\dx\,$-closed form $\alpha\in \Omega^{1,1}_{\Rnum}$, representing an element of $\tilde{H}^{1,1}_{\Rnum}$. Then $0= \dx \alpha = \del \alpha + \delb \alpha$. 
    Since the two summands are of different bidegree, they must both be zero, and so as a complex form $\alpha$ is also  $\del$- and $\delb$-closed. 

    If $\alpha = \dx \eta$ for some $\eta\in \Omega^1(M,\Rnum)$, by the $\del\delb$-Lemma $\alpha$ is also $\delb$-exact. 

    Conversely, a $\delb$-closed real $\alpha\in \Omega^{1,1}_{\Rnum}$ is also $\dx\,$-closed: $\dx \alpha= \del \alpha + \delb \alpha= \overline{\delb \alpha} + \delb \alpha =0$. The $\del\delb$-Lemma then again implies that the notions of $\delb$- and $\dx\,$-exactness also agree. 
\item 
For a real form $\beta = \beta^{2,0}+\beta^{0,2} \in \left(\Omega^{2,0}\oplus \Omega^{0,2}\right)_{\Rnum}$, by definition: $\beta^{0,2}=\overline{\beta^{2,0}}$.

We have  
    \[
       \dx \beta =\del\beta^{2,0} + \delb\beta^{2,0} + \del \beta^{0,2}  + \delb \beta^{0,2}. 
    \] Again, since the summands on the right-hand side all have different complex bidegree, the following are equivalent: 
    \begin{enumerate}[label=(\roman*)]
        \item $\beta$ is $\dx\,$-closed. 
        \item $\beta^{2,0}, \beta^{0,2}$ are $\dx\,$-closed.
        \item $\beta^{2,0}, \beta^{0,2}$ are $\delb$-closed. 
    \end{enumerate}

    Now consider $\lambda\in\Omega^1(M,\Rnum)$. Split into distinct $(p,q)$-types, we have $\lambda=\lambda^{1,0}+\lambda^{0,1}$ with $\lambda^{0,1}=\overline{\lambda^{1,0}}$.  
    Evaluate: 
     \begin{equation*}
        \dx \lambda = \del\lambda^{1,0}  +\delb \lambda^{1,0} + \del \lambda^{0,1} + \delb \lambda^{0,1}. 
    \end{equation*}

If $\dx \lambda \in \left(\Omega^{2,0}\oplus \Omega^{0,2}\right)_{\Rnum}$, in order to match $(p,q)$-type, we must have: 
\[ \dx \lambda = \del\lambda^{1,0} + \delb \lambda^{0,1}. \]

So a 
%$\dx\,$-closed 
form $\beta$ in $\left(\Omega^{2,0}\oplus \Omega^{0,2}\right)_{\RR}$ is $\dx\,$-exact if and only if its $(2,0)$-component $\beta^{2,0}$ is $\del$-exact. 
By the $\del\delb$-Lemma applied to $\beta^{2,0}$, this is then the case if and only if $\beta^{2,0}$ is also $\dx\,$- and $\delb$-exact. In that case, 
$\beta^{0,2}$ is also $\delb$-exact. (Since its bidegree is $(2,0)$, that actually means $\beta^{2,0}=0$.)

This shows that the map 
\[ 
\tilde{H}^{(2,0)+(0,2)}_{\Rnum} \to
\left(H^{2,0}\oplus H^{0,2}\right)_{\Rnum},
\quad
[\beta]\mapsto [\beta^{2,0}]+[\beta^{0,2}]\]
is well-defined and injective. The first statement made in this item assures that the map is surjective, and hence gives the canonical isomorphism required by the lemma. 
\end{enumerate}
    
We have thus shown that for both cases of definite real bidegree, the suitably-defined $\delb$- and $\dx\,$-cohomologies agree. 
\end{proof}

Using the orthogonality at the level of forms, Lemma \ref{lem:HRdclosed} and the canonical decomposition of complex de Rham into Dolbeault cohomology given in \eqref{eq:Hodge} it thus follows that: 

\begin{cor}\label{cor:realsplitting}
    The splitting of de Rham cohohomology 
    \[H^2(M,\Rnum)=\tilde{H}^{(2,0)+(0,2)}_{\Rnum}  \oplus \tilde{H}^{(1,1)}_{\Rnum}\]
    is orthogonal with respect to the inner product on $H^2(M,\Rnum)$ given by integration. 
\end{cor}

\appendix

\bibliographystyle{habbrv} 
%\bibliography{brane}

\end{document}